\documentclass[10pt,a4paper]{amsart}
\usepackage{amsmath,amssymb}
\usepackage[latin1]{inputenc}
\usepackage[english]{babel}

\let\8=\infty
\let\landa=\lambda

\let\parc=\partial

\def\landa{\lambda}

\def\flecha{\rightarrow}
\def\esiz{\langle}
\def\esde{\rangle}

\newcommand\esc[2]{\langle{#1},{#2}\rangle}

\def\cosh{\mathop{\rm cosh }\nolimits}

\def\S{\mathbb{S}}

\def\E{\mathbb{E}}

\def\R{\mathbb{R}}
\def\C{\mathbb{C}}
\def\H{\mathbb{H}}
\def\X{\mathfrak{X}}

\def\Ek{\mathbb{E}^3 (\kappa,\tau)}

 \newtheorem{defi}{Definition}[section]
 \newtheorem{teo}[defi]{Theorem}
 \newtheorem{pro}[defi]{Proposition}
 \newtheorem{cor}[defi]{Corollary}
 
 \newtheorem{eje}[defi]{Example}
 \newtheorem{remark}[defi]{Remark}

\numberwithin{equation}{section}

\begin{document}

\title[CMC surfaces in homogeneous $3$-manifolds]{A characterization of constant mean curvature surfaces in homogeneous $3$-manifolds}

\author{Isabel Fernández}

\address{Departamento de Geometría y Topología, Universidad de Granada,
E-18071 Granada, Spain} \email{isafer@ugr.es}

\author{Pablo Mira}
\address{Departamento de Matemática Aplicada y Estad\' \i stica, Universidad Politécnica
de Cartagena, E-30203 Cartagena, Murcia, Spain} \email{pablo.mira@upct.es}


\date{}


\begin{abstract}
It has been recently shown by Abresch and Rosenberg that a certain Hopf differential is holomorphic on every constant mean curvature surface in a Riemannian homogeneous $3$-manifold with isometry group of dimension $4$. In this paper we describe all the surfaces with holomorphic Hopf differential in the homogeneous $3$-manifolds isometric to $\H^2\times \R$ or having isometry group isomorphic either to the one of the universal cover of ${\rm PSL} (2,\R)$, or to the one of a certain class of Berger spheres. It turns out that, except for the case of these Berger spheres, there exist some exceptional surfaces with holomorphic Hopf differential and non-constant mean curvature.
\end{abstract}

\maketitle

\section{Introduction}

An extremely useful tool in surface theory is the fact that the Hopf differential of a surface in a $3$-dimensional space form is holomorphic if and only if the surface has constant mean curvature (CMC). Inspired by this result, Abresch and Rosenberg proved in \cite{AbRo} that for CMC surfaces in the product spaces $\H^2\times \R$ and $\S^2\times \R$ there is a certain \emph{perturbed} Hopf differential which is holomorphic. This differential may be seen as the usual Hopf differential of the surface plus a certain correction term. Even more generally, Abresch showed in \cite{Abr} the existence of such a holomorphic Hopf-type differential for CMC surfaces immersed in $3$-dimensional homogeneous manifolds with $4$-dimensional isometry group. These results have made of CMC surfaces in homogeneous $3$-manifolds a fashion research topic, on which many interesting works are being produced at the present time. An \emph{almost} up-to-date reference list on this subject may be consulted in \cite{FeMi}.

There is, however, a natural question that remains unanswered: \emph{are CMC surfaces in homogeneous $3$-manifolds the only surfaces for which the Hopf-type differential introduced by Abresch and Rosenberg is holomorphic?} In other words, one wishes to know if the converse of the above mentioned results by Abresch and Rosenberg hold. This has been a frequently discussed problem among people working on this area. The only known particular solution to this problem was found by Berdinsky and Taimanov in \cite{BeTa}, where is it proved that the converse holds when the homogeneous target space is the $3$-dimensional Heisenberg group ${\rm Nil_3}$.

In this paper we will give an answer to the above question for certain homogeneous $3$-manifolds. But first, in order to state our result, some basic comments on Riemannian homogeneous $3$-manifolds should be made. The details may be consulted in \cite{Abr,BeTa,Dan}, for instance.

The homogeneous $3$-manifolds with $4$-dimensional isometry group can be classified in terms of a pair of real numbers $(\kappa,\tau)$ satisfying $\kappa\neq 4\tau^2$. Indeed, all these manifolds are fibrations over a complete simply-connected surface $\mathcal{M}^2(\kappa)$ of constant curvature $\kappa.$ Translations along the fibers are isometries and therefore they generate a Killing field, $\xi,$ also called  {\em the vertical field}. The number $\tau$ is the one such that 
$ \overline{\nabla}_X \xi=\tau X\times\xi $
holds for any vector field $X$ on the manifold. Here $\overline{\nabla}$ is the Levi-Civita connection of the manifold and $\times$ denotes the cross product.

It is important to notice that for $\tau=0$ this fibration becomes trivial and thus we get the product spaces $\mathcal{M}^2(\kappa)\times\mathbb{R}$. When $\tau\neq 0$ the manifolds have the isometry group of the Heisenberg space if $\kappa=0,$ of the Berger spheres if $\kappa>0$, or the one of the universal covering of $\mbox{PSL}(2,\mathbb{R})$ when $\kappa<0$.

In what follows $\mathbb{E}^3(\kappa,\tau)$ will represent a homogeneous 3-manifold with isometry group of dimension 4, where $\kappa$ and $\tau$ are the real numbers described above.

For an immersion $\psi:\Sigma\to\Ek$, let $pdz^2$ be its Hopf differential, i.e. the $(2,0)$ part of its complexified second fundamental form. Then there exists a quadratic differential $Ldz^2$ defined in terms of $\kappa,$ $\tau$, the mean curvature $H$ of $\psi$ and the restriction of the vertical field $\xi$ on the surface, such that $Qdz^2:=pdz^2+Ldz^2$ is holomorphic whenever $H$ is constant \cite{Abr,AbRo} (see Section 2 for the details). We shall call $Qdz^2$ the {\em Abresch-Rosenberg differential} of the surface.

With all of this, our main result is the following:

\begin{teo}\label{main}
Let $\Ek$ be a homogeneous $3$-manifold of base curvature $\kappa$ and bundle curvature $\tau$, and let $\psi:\Sigma\flecha \Ek$ be a surface with holomorphic Abresch-Rosenberg differential. Then
 \begin{enumerate}
 \item
If $0<\kappa /8 \leq \tau^2$, i.e. the group of isometries of $\Ek$ is isomorphic to the one of a Berger sphere of a certain type, then $\psi$ is a CMC surface. 
 \item
If $\kappa <0$ and $\tau =0$ (i.e. $\Ek \equiv \H^2 (k)\times \R$), then $\psi$ is a CMC surface, or it is one of the rotational surfaces in Example \ref{ex2}.
 \item
If $\kappa <0$ and $\tau <0$ (i.e. $\Ek$ has isometry group isomorphic to the one of the universal cover of ${\rm PSL} (2,\R)$), then $\psi$ is a CMC surface, or it is one of the surfaces in Example \ref{ex3}.
 \end{enumerate}
\end{teo}

Some remarks should be made regarding our result.
 \begin{enumerate}
 \item
The proof of this theorem will also show that the surfaces with holomorphic Abresch-Rosenberg differential in the Heisenberg $3$-space are CMC surfaces. This was proved in \cite{BeTa}.
 \item
It remains unsolved whether CMC surfaces are the only surfaces in $\S^2\times \R$ or in the Berger spheres $\Ek$ with $0<8 \tau^2 < \kappa$ that have holomorphic Abresch-Rosenberg differential.
 \item
We shall prove that a compact surface in $\Ek$ with holomorphic Abresch-Rosenberg differential (and non-zero Euler characteristic if $\tau \neq 0$) is always a CMC surface.
 \end{enumerate}

The outline of the paper goes as follows. In Section 2 we will describe the integrability equations for surfaces in the homogeneous $3$-manifolds $\Ek$ in terms of an isothermal coordinate patch. We hope that this approach will be useful for the application of integrable systems techniques to the study of CMC surfaces in these spaces. We shall also show in Section 2 that a surface with vanishing Abresch-Rosenberg differential is a CMC surface.
In Section 3 we will expose some exceptional surfaces in certain homogeneous spaces $\Ek$ that have holomorphic Abresch-Rosenberg differential, but which have non-constant mean curvature. In Section 4 we shall prove Theorem \ref{main}. We will also show there that a compact surface with holomorphic Abresch-Rosenberg differential (and non-zero Euler characteristic if $\tau \neq 0$) is always a CMC surface.

We finally wish to point out that the techniques of this paper can be used in some other related geometrical theories. For instance, a holomorphic quadratic differential for surfaces of constant curvature in $\H^2\times \R$ and $\S^2\times \R$ has been recently found in \cite{AEG}, and it is likely that our considerations here may be extended to that setting.

\section{Surfaces in homogeneous $3$-manifolds}

In this section we will describe the fundamental equations for an
immersed surface $\psi:\Sigma\to\Ek$ in terms of a conformal
parameter $z$ on the surface. So, we will consider $\Sigma$ as a Riemann surface with the conformal structure given by its induced metric via $\psi$, and we will let $z$ denote a conformal parameter of $\Sigma$. Associated to $z=s+it$, we will consider the usual operators $\parc_z = (\parc_s -i\parc_t )/2$ and $\parc_{\bar{z}} =(\parc_s +i \parc_t)/2$. With all of this, we will define the following {\em fundamental data}.

\begin{defi}\label{fundada} In the above setting, let $\eta$
be the unit normal map of $\psi$, and let $\xi$ denote the vertical unit Killing field of $\E^3(\kappa,\tau)$. We will call the \emph{fundamental
data} of $\psi$ to the uple $(\lambda,u,H,p,A)\in\R^+\times
[-1,1]\times\R\times\C^2$, where

\begin{tabular}{l}
\qquad$\lambda$  is the conformal factor of the induced metric in
$\Sigma,$
$\lambda=2\esc{\psi_z}{\psi_{\bar z}}.$\\
\qquad$u$  is the normal component of the vertical field $\xi$,
$u=\esc{\eta}{\xi}.$\\
\qquad$H$  is the mean curvature of $\psi$.\\
\qquad$p\, dz^2$ is the Hopf differential of $\psi$, $p=-\esc{\psi_z}{\eta_z}$.\\
\qquad$A$  $= \esc{\xi}{\psi_z}=\esc{T}{\partial_z}$ where $T\in\mathfrak{X}(\Sigma)$ is given by
$d\psi(T)=\xi-u\eta.$
\end{tabular}
\end{defi}

\begin{remark}\label{re:A} When $\tau=0$ (i.e., $\Ek=\mathcal{M}^2(\kappa)\times\R$) the vertical field $\xi$ is nothing but 
$\xi=(0,1)\in\mathfrak{X}(\mathcal{M}^2(\kappa))\times\R$. Therefore, if we write $\psi=(N,h):\Sigma\to\mathcal{M}^2(\kappa)\times\R$ then $A=h_z$.
As a consequence, if $A$ is identically zero the surface is a piece of a horizontal slice, which has $H=0$. However, if $\tau\neq 0$ $A$ cannot vanish on an open subset of $\Sigma$ (see equations {\bf (C.2)} and {\bf (C.4)} in Theorem \ref{th:formulas}).
\end{remark}

\begin{teo}\label{th:formulas}
The fundamental data of an immersed surface $\psi:\Sigma\flecha \Ek$ satisfy the following integrability conditions:

\begin{equation}\label{lasces}
\left\{\def\arraystretch{1.3} \begin{array}{lccc} {\bf (C.1)} & p_{\bar{z}} & = 
& \displaystyle \frac{\landa}{2} (H_z + u A (\kappa - 4\tau^2)). \\ {\bf (C.2)} & 
A_{\bar{z}} & = & \displaystyle \frac{u \landa}{2} (H+i\tau) .\\ {\bf (C.3)} & 
u_{z} & = & - (H-i\tau) A -\displaystyle \frac{2 p}{\landa} \bar{A}.\\ {\bf 
(C.4)} & \displaystyle \frac{4 |A|^2}{\landa} & = & 1 - u^2 .
\end{array}\right. 
\end{equation} 

Conversely, let us choose functions $\landa,u,H :\Sigma\flecha \R$ with $\landa >0,$ $-1\leq u \leq 1$ and $p,A:\Sigma\flecha \C$ on a simply connected Riemann surface $\Sigma$, verifying \eqref{lasces} for some real constants $\kappa,\tau$ with $\kappa - 4\tau^2 \neq 0$.

Then there exists a unique (up to congruences) surface $\psi :\Sigma\flecha \Ek$ whose fundamental data are $\{\landa,u,H,p,A\}$.
\end{teo}

\begin{proof} The proof of this theorem follows from Theorem 4.3 in
\cite{Dan} where a necessary and sufficient condition for the existence of an isometric immersion from a Riemannian surface $\Sigma$ into $\mathbb{E}^3(\kappa,\tau)$ is given
in terms of the following {\em compatibility equations} for any
vector fields $X,$ $Y$ on $\Sigma$:

\begin{equation*}
\begin{array}{lc}
(i)& \|T\|^2+u^2=1,\\[0,2cm]
(ii)& du(X)+\esc{SX-\tau JX}{T}=0,\\[0,2cm]
(iii)& \nabla_X T=u(SX-\tau JX),\\[0,2cm]
(iv)& \nabla_X SY- \nabla_Y SX
-S[X,Y]=(\kappa-4\tau^2)u\big(\esc{Y}{T}X-\esc{X}{T}Y\big),\\[0,2cm]
(v)& K=\mbox{det}S+\tau^2+(\kappa-4\tau^2)u^2.
\end{array}
\end{equation*}

Here $\nabla$ denotes the Levi-Civita connection on the surface,
$K$ is its Gauss curvature, $S$ is the shape operator and $J$ is
the complex structure on $\Sigma$. Thus, we just need to check
that our fundamental equations
{\bf (C.1)} to {\bf (C.4)} are equivalent to the above ones.

More specifically, let $\psi:\Sigma\flecha \Ek$ be a surface with fundamental data $(\landa,u,H,p,A)$ with respect to a conformal parameter $z$ of $\Sigma$. We will keep the above notations for $T,S,\nabla,K$ and $J$. Then, by the definition of $p$ and $H$ we have 
 \begin{equation}\label{sure1}
 \esiz S(\parc_z),\parc_z\esde = p,\hspace{1cm} \esiz S(\parc_z),\parc_{\bar{z}}\esde = \frac{\landa H}{2}.
 \end{equation}
Moreover, since $A=\esiz T,\parc_z\esde$, we can write 
 \begin{equation}\label{sure2}
 T=\frac{2}{\landa} \left(\bar{A} dz + Ad \bar{z} \right).
 \end{equation}
In addition, we also have the following intrinsic metric relations on $\Sigma$:
\begin{equation}\label{sure3}
\def\arraystretch{1.7} \begin{array}{lll}
\esiz \parc_z,\parc_z\esde =0 & \esiz \parc_z,\parc_{\bar{z}}\esde =\displaystyle\frac{\landa}{2} & K=\displaystyle\frac{-2 (\log \landa )_{z\bar{z}}}{\landa} \\ \nabla_{\parc_z} \parc_z = \displaystyle \frac{\landa_z}{\landa} \parc_z & \nabla_{\parc_z} \parc_{\bar{z}} = 0 & J(\parc_z)= i \parc_z.
 \end{array}
\end{equation} 

With this, let $(\landa,u,H,p,A)$ be an uple in the conditions of Definition \ref{fundada}, and let us consider in terms of them the Riemannian surface $(\Sigma,\landa |dz|^2)$, whose fundamental data are given by \eqref{sure3}, the symmetric endomorphism $S:\X(\Sigma)\flecha \X (\Sigma)$ described by \eqref{sure1} and the unit tangent field $T\in \X(\Sigma)$ of \eqref{sure2}. We are going to show that $(\landa,u,H,p,A)$ satisfy {\bf (C.1)} to {\bf (C.4)} if and only if $(S,T,\nabla,u,J,K)$ verify $(i)$ to $(v)$. This will finish the proof, by Theorem 4.3 in \cite{Dan}.

First of all, it is direct to observe that, by \eqref{sure2}, 
$(i)$ is equivalent to ${\bf (C.4)}$.

Secondly, $(ii)$ is a linear expression in the variable $X$,
and hence it suffices to show that this equality holds for
$X=\partial_z.$ But this coincides with ${\bf (C.3)}$ so both
equations are also equivalent.

Likewise, we just need to check that $(iii)$ holds for
$X=\partial_z$. This is equivalent to show that
\begin{equation}\label{eq:esc}
\esc{\nabla_{\partial_{z}}T}{\partial_{\bar
z}}=\frac{u\lambda}{2}(H-i\,\tau) \quad \mbox{and}\quad
\esc{\nabla_{\partial_z}T}{\partial_z}=u\,p.
\end{equation} 
By using the fact that $ \esiz \nabla_X T,Y\esde =X\esiz T,Y\esde - \esiz T,\nabla_X Y\esde$, as well as the identities in \eqref{sure3},  it is not hard to see that the
first equation in \eqref{eq:esc} is precisely the conjugate expression of {\bf (C.2)}, while the second one can
be rewritten as
$${\bf (C.0)} \qquad\quad A_z-\frac{\lambda_z}{\lambda}A = up.$$
This equation {\bf (C.0)} is obtained by deriving ${\bf (C.4)}$
with respect to $z$ and using ${\bf (C.2)}$ and ${\bf
(C.3)}$.


Our aim now is to show that $(iv)$ is equivalent to {\bf (C.1)}.
Indeed, $(iv)$ is an anti-symmetric bilinear expression in the
variables $X$,$Y$ and so it suffices to check it for $X=\partial_z$,
$Y=\partial_{\bar z}$. Moreover, it is enough to see that if we take
scalar product with $\partial_z$ then the equality holds (by taking
scalar products with $\partial_{\bar z}$ we get the conjugate
expression). Using the same ideas as above, and after some calculations, we see that 
this is exactly {\bf (C.1)}, so we are done.

Finally, taking into account that
$$ K=\frac{-2(\log\lambda)_{z\bar z}}{\lambda} \qquad \mbox{and}
\qquad \mbox{det}S=H^2-\frac{4|p|^2}{\lambda^2}$$ we can write
$(v)$ as $$ (\mbox{log}\,\lambda)_{z\bar
z}=\frac{2|p|^2}{\lambda}-\frac{\lambda}{2} u^2
(\kappa-4\tau^2)-\frac{\lambda}{2}(H^2+\tau^2).$$ Straightforward
computations show that this equation is obtained by deriving {\bf
(C.0)} with respect to $\bar z$ and using {\bf (C.1)}, {\bf (C.2)} and {\bf
(C.3)}. This finishes the proof.
\end{proof}


The Codazzi equation {\bf (C.1)} can be expressed in an alternative way that will be more useful to us. For this, let us first define for an immersed surface $\psi:\Sigma\flecha \Ek$ its \emph{Abresch-Rosenberg} differential as the quadratic differential $$Q dz^2 = \left( 2 p - \frac{\kappa -4 \tau^2}{H+i\tau } A^2 \right) dz^2,$$ following the above notations, and defined away from points with $H=0$ if $\tau =0$. We will assume from now on that the surfaces have non-vanishing mean curvature if $\tau=0$. There is no loss of generality with that, since our study is basically local. We are just excluding the minimal surfaces in $\S^2\times \R$ and $\H^2\times \R$, that are better studied by other methods.

It is then easy to see by means of ${\bf (C.2)}$ that the Codazzi equation can be rephrased in terms of $Q$ as 
 \begin{equation}\label{cod}
  Q_{\bar{z}}= \landa H_z + (\kappa - 4\tau^2 ) \frac{H_{\bar{z}} A^2}{(H+i\tau)^2}.
  \end{equation}
Consequently, one has
 \begin{cor}\cite{Abr,AbRo}\label{abro} $Q dz^2$ is a holomorphic quadratic differential on any CMC surface in $\Ek$.
  \end{cor}

Our purpose in this work is to describe to what extent the holomorphicity of the Abresch-Rosenberg differential characterizes the CMC surfaces in the homogeneous $3$-manifolds $\Ek$. As a preliminary step for this, let us describe first of all the case in which $Q$ vanishes identically. 

%

\begin{pro}\label{lema}
Any surface in $\Ek$ with vanishing Abresch-Rosenberg differential is a CMC surface.
\end{pro}
\begin{proof}
Assume that $H$ is non-constant in some open set of $\Sigma$. Then we may suppose without loss of generality that $H_z\neq 0$ and $A\not\equiv 0$. As $Q\equiv 0$, by its own definition we obtain 
 \begin{equation}\label{1lema}
 2p =\frac{\kappa - 4\tau^2 }{H+i \tau} A^2.
 \end{equation}
Taking modulus on the Codazzi equation \eqref{cod}, and using $H_z\neq 0$ we find that 
 \begin{equation}\label{3lema}
 \frac{|A|^2}{\landa} = \frac{H^2 +\tau^2}{|\kappa- 4\tau^2|}.
 \end{equation}
If we substitute now equations \eqref{1lema} and \eqref{3lema} into ${\bf (C.3)}$ we end up with 
 \begin{equation}\label{35lema}
 - u_z =(H-i\tau) A \left( 1 +\frac{\kappa -4\tau^2}{|\kappa - 4\tau^2|}\right).
  \end{equation} Consequently, if $\kappa - 4\tau^2 <0$ we infer that $u$ is constant. But on the other hand, putting together \eqref{3lema} and ${\bf (C.4)}$ we obtain 
 \begin{equation}\label{4lema}
 1-u^2 = \frac{4 (H^2 +\tau^2)}{|\kappa - 4\tau^2|}.
 \end{equation}
So $H$ should be constant in this case, and this is not possible. Hence $\kappa-4\tau^2 >0$. This indicates via \eqref{cod} that $H_z (H+i\tau) \bar{A} \in i\R$. Therefore, by differentiating \eqref{4lema} we find that $u_z (H+i\tau)\bar{A}\in i\R$. But this is not possible, by \eqref{35lema}. This concludes the proof.

\end{proof}


\begin{remark}
The way we have defined the quadratic differential $Q dz^2$ is not exactly the way it was defined in \cite{AbRo} and \cite{Abr}. Indeed, in these works the authors work with the differential $P dz^2 = (H+i\tau) Q dz^2$. This obviously makes no difference when working with CMC surfaces, but it does in the present situation. The definition of the Abresch-Rosenberg differential we have adopted here is taken from Berdinsky-Taimanov \cite{BeTa}, which is the first paper in where the converse of Corollary \ref{abro} is treated. Indeed, for the case $\tau \neq 0$, the definition of the Abresch-Rosenberg differential as $Q dz^2$ is from a certain viewpoint more natural, as it is constructed by adding to the usual Hopf differential $p dz^2$ a certain \emph{correction term} $L dz^2$. 

We will discuss briefly the situation for the quadratic differential $P dz^2$ in the product spaces $\S^2\times \R$ and $\H^2\times \R$ in Example \ref{lacompli}.
\end{remark}


\section{Some exceptional examples}

\begin{eje}\label{lacompli}
Let $a,b\in\mathbb{R}$ be two real constants with $a\neq 0,$ set $\kappa=\pm 1,$ and consider the real function $h(s)$ given by 
\begin{equation}\label{eq:h}
h'(s)=\left\{\def\arraystretch{1.8}\begin{array}{ll}
\displaystyle{\frac{-1}{\sinh(as+b)}} & \mbox{if}\; \kappa=-1,\\
\displaystyle{\frac{1}{\cosh(as+b)}} & \mbox{if}\; \kappa=1,
\end{array}\right.
\end{equation}
defined for $as +b>0$. So, $h'(s)$ is the general solution of the autonomous ODE
\begin{equation}\label{eq:eqdifh}
y'=- a\,y\,\sqrt{1-\kappa\, y^2}.
\end{equation}
Next, define in terms of $h'(s)$ the following quantities
\begin{equation}\label{eq:def}\begin{array}{cc}
u=\displaystyle{\frac{a}{\sqrt{1+a^2}}}, & \qquad\lambda=(1+a^2)(h')^2,\\[0,5cm]
H=-\displaystyle{\frac{1}{2\sqrt{1+a^2}}}\sqrt{-\kappa+\frac{1}{(h')^2}}, &\qquad p=\displaystyle{\frac{-\lambda\,H}{2}}.
\end{array}\end{equation}

It follows by an elementary computation using (\ref{eq:eqdifh}) and (\ref{eq:def}) that the quantities $\{ \lambda, u, H, p, A:=h'/2\}$ verify equations {\bf (C.1)} to {\bf (C.4)} for $\tau=0$ in terms of the complex parameter $z=s+it$, where $t$ is an arbitrary real parameter. As a consequence, we obtain a surface immersed in $\mathbb{E}^3(\kappa,0)=\mathcal{M}^2(\kappa)\times\mathbb{R}$ having non-constant mean curvature $H$. Moreover, 
$$ Pdz^2=HQdz^2=\frac{-1}{4}dz^2$$
is holomorphic, even though $H$ is non constant.\\

At last, let us point out that the surfaces described in this way are rotational in $\mathcal{M}^2(\kappa)\times\mathbb{R}.$
For any $t_0\in\mathbb{R}$ the map $(s,t)\mapsto (s,t+t_0)$ preserve all the fundamental data $\{\lambda(s),u(s),p(s),H(s),A(s)\}$
and consequently $\psi$ has a continuous $1$-parameter group of self-congruences. In other words, for every $t_0\in\mathbb{R}$ there is a rigid motion $\Psi_{t_0}$ of $\mathcal{M}^2(\kappa)\times\R$ satisfying $\psi(s,t+t_0)=\Psi_{t_0}(\psi(s,t))$. So, $\{ \Psi_{t_0}\;:\;t_0\in\mathbb{R}\}$ is a continuous $1$-parameter group of isometries of $\mathcal{M}^2(\kappa)\times\R$ i.e., it consists of
helicoidal motions. At last, as $h(s,t+t_0)=h(s,t)$ and $h$ is precisely the last coordinate function of the immersion $\psi$ (see Remark \ref{re:A}), we conclude that all the $\Psi_{t_0}$'s are rotations. Thus the surface is rotationally invariant.

\end{eje}


\begin{eje}\label{ex2}
Let $\alpha=\alpha(r)$ be a solution of the following second order autonomous ODE:
\begin{equation}\label{eq:3.2} 
\alpha''=(\alpha')^2\,\cot \alpha -\delta \cos \alpha, 
\end{equation}
where $\delta=\pm 1.$ Define the following data:
\begin{equation}\label{def3.2}\def\arraystretch{3}\begin{array}{lll}
\lambda =\displaystyle\frac{1}{(\alpha')^2} ,&  H=\displaystyle\frac{\sin\alpha}{2}, &
u=\cos\alpha, \\ 
p=\displaystyle\frac{1}{2}-\displaystyle\frac{\delta\sin\alpha}{4(\alpha')^2} ,&
\hspace{0.5cm} A=\displaystyle \frac{\delta^{1/2}\sin\alpha}{2\alpha'} , &
\end{array} 
\end{equation}
whenever $\alpha'(r)\neq 0$.
Straightforward computations show that these data satisfy equations {\bf (C.1)} to {\bf (C.4)} with $\kappa=-1,$ $\tau=0,$ for the 
parameter $z=r+i\,\theta$ if $\delta=1,$ or $z=\theta+i\,r$ if $\delta=-1.$ Here $\theta$ is a real parameter.
Consequently, by Theorem 3.1 they are the fundamental data of an immersed surface in $\mathbb{H}^2\times\mathbb{R}$ with 
$Q\,dz^2=dz^2$ and non-constant mean curvature. Moreover, using the argument in Example 3.1, and bearing in mind that $A\in \R$ if $\delta =1$ and $A\in i\R$ if $\delta =-1$, we infer that this surface is rotationally invariant.

\end{eje}


\begin{eje}\label{ex3}
Let $\kappa,\tau\in\R$ with $\kappa<0,$ $\tau\neq 0$ and $\kappa-4\tau^2<0.$ Let us consider $z=s+it$ a complex parameter and $H(s,t)$ a non-constant solution of the following overdetermined system of PDEs:
\begin{equation}\label{eq:Hdif}
\left\{\begin{array}{l}
\big(\log (H^2+\tau^2)\big)_{z\bar z}= \displaystyle{\frac{8H^2|H_z|^2}{(H^2+\tau^2)(4H^2+\kappa)}}+
\displaystyle{\frac{H_z^2(H+i\tau)(4H^2+\kappa)}{4|H_z|^2(H^2+\tau^2)}},\\[0,5cm]
H_z^2(H+i\tau) \in\R \qquad \mbox{i.e.,} \quad \tau(H_s^2-H_t^2)=2HH_s H_t,
\end{array}\right.
\end{equation}
with the condition $4H^2+\kappa<0$ and $H_z\neq 0$.

Define next in terms of $H$ the quantities:
\begin{equation}\label{eq:def3}
\begin{array}{cc}
u=\sqrt{\displaystyle{\frac{4H^2+\kappa}{\kappa-4\tau^2}}}, &\qquad A=\displaystyle{\frac{u(H^2+\tau^2)}{4HH_z}},\\[0,5cm]
p=\displaystyle{\frac{1}{2}\big(1+\frac{\kappa-4\tau^2}{H+i\tau} A^2\big)}, &\qquad \lambda= \displaystyle{\frac{-|A|^2(\kappa-4\tau^2)}{H^2+\tau^2}}.
\end{array}
\end{equation}

Then, an elementary computation indicates that $\{\lambda, u, H, p, A\}$ satisfy conditions {\bf (C.1)}, {\bf (C.3)} and {\bf (C.4)} in Theorem \ref{th:formulas}. Moreover, they also verify {\bf (C.2)}. To see this, we first observe that by \eqref{eq:def3} $$\frac{u}{2A} = \frac{2 H H_z}{H^2 +\tau^2} =(\log (H^2 +\tau^2))_z.$$ With this, by deriving this expression with respect to $\bar{z}$, and using {\bf (C.3)} together with the first formula in \eqref{eq:Hdif}, we obtain {\bf (C.2)}.

As a consequence, the data in \eqref{eq:def3}, \eqref{eq:Hdif} define a surface in $\Ek$ with holomorphic Abresch-Rosenberg differential, $Qdz^2=dz^2$, and non-constant mean curvature $H$. To finish this example, we need to ensure that the system \eqref{eq:Hdif} has some solution. In order to do so, let us define $f(x)=-\frac{x}{\tau}\pm\sqrt{1+(\frac{x}{\tau})^2}$. Then the second equation in \eqref{eq:Hdif} can be written as 
\begin{equation}\label{eq:Ht}
H_t=f(H)H_s.
\end{equation}
This equation has two basic consequences. On one hand, by inserting this relation into the first equation in \eqref{eq:Hdif}, and after some computations, we find that
\begin{equation}\label{eq:Hss}
H_{ss}=\mathcal{F}(H,H_s,\kappa,\tau)
\end{equation}
for some real analytic function $\mathcal{F}$. On the other hand, (\ref{eq:Ht}) can be solved, and the general solution is given by the implicit relation 
\begin{equation}\label{eq:implicit}
s+t\,f(H)=g(H),
\end{equation}
where $g$ is an arbitrary smooth real function.
By differentiating this relation we get, 
\begin{equation}\label{eq:7}\displaystyle{
H_s=\frac{1}{g'(H)-tf'(H)}\quad\mbox{and}\quad  H_{ss}=-\frac{g''(H)-t\,f''(H)}{(g'(H)-t\,f'(H))^3} }.
\end{equation}
So, plugging these expressions into \eqref{eq:Hss} it is obtained 
\begin{equation}\label{eq:gdif}
g''(H) = t\,f''(H)-\big(g'(H)-tf'(H)\big)^3\mathcal{F}\big(H,\frac{1}{g'(H)-t\,f'(H)},\kappa,\tau\big).
\end{equation}
In other words, $g(x)$ is a solution of an ODE of the form 
\begin{equation}\label{eq:gdif2}
g''=\mathcal{G}(H,g,g',\kappa, \tau,t), 
\end{equation}
where here $t$ is considered as a real parameter, and $\mathcal{G}$ is an analytic function induced by $\mathcal{F}$ via 
\eqref{eq:gdif}.

Now, given $g(x)=g(x,\kappa, \tau, t)$ a solution of \eqref{eq:gdif2}, consider $H(s,t)$ given implicitly by \eqref{eq:implicit} (this can be always achieved locally if we assume that $g'(x)\neq 0$). Then by \eqref{eq:7} and \eqref{eq:gdif} together with \eqref{eq:implicit} we find that $H$ verifies \eqref{eq:Ht} and \eqref{eq:Hss}. Consequently, $H$ must satisfy \eqref{eq:Hdif}. This ensures the existence of these examples.

\end{eje}


\section{Proof of Theorem \ref{main}}

Let $\psi:\Sigma\flecha \Ek$ be an immersed surface with holomorphic Abresch-Rosenberg differential $Q dz^2$, and let us assume that it is not a CMC surface. By Lemma \ref{lema}, $Q$ does not vanish identically, and consequently $Q$ has isolated zeros. So, working in a simply connected piece of $\Sigma$ away from these zeros, it is possible to introduce a new complex parameter (which will also be denoted by $z$) so that $Q dz^2 \equiv dz^2$, i.e. $Q\equiv 1$. Therefore it holds
 \begin{equation}\label{main1}
 2p-1 =\frac{\kappa - 4\tau^2}{H+i\tau} A^2.
 \end{equation}
Besides, the Codazzi equation \eqref{cod} gives 
 \begin{equation}\label{main3}
 -\landa |H_z|^2 = (\kappa -4\tau^2)\left(\frac{A H_{\bar{z}}}{H+i\tau}\right)^2.
 \end{equation}
Taking modulus in (\ref{main3}) and using {\bf (C.4)} as well as $H_z\neq 0$, we get 
\begin{equation}\label{main4}
1\geq 1-u^2= \frac{4|A|^2}{\landa} =  \frac{4 (H^2 +\tau^2)}{|\kappa -4\tau^2|}.
\end{equation}
This inequality implies that if $\kappa \geq 0$ and $\kappa -4\tau^2 <0$, the surface $\psi$ cannot exist. The same can be said if $\kappa -4\tau^2 >0$ but $\kappa -8\tau^2 \leq 0$. In this way we have proved the first part of the theorem. We have also shown at this point that all surfaces with holomorphic Abresch-Rosenberg differential in the Heisenberg $3$-space ($\kappa =0$) are CMC surfaces (this is a result in \cite{BeTa}).

Let us assume from now on that $\kappa - 4\tau^2 <0$. Then \eqref{main4} provides $4H^2 + \kappa \leq 0$ and 
 \begin{equation}\label{main5}
 u=\sqrt{\frac{4H^2+\kappa}{\kappa -4\tau^2}}
 \end{equation}
(up to a $\pm$ sign, that can be changed by reversing orientation if necessary). Moreover, differentiation of \eqref{main4} gives us 
 \begin{equation}\label{main6}
4H H_z = (\kappa -4\tau^2) u u_z.
 \end{equation}
Now, putting together ${\bf (C.3)}$, \eqref{main1} and \eqref{main4}, we obtain 
 \begin{equation}\label{main6p5}
u_z = -\bar{A} /\landa, 
 \end{equation}
and so \eqref{main4} and \eqref{main6} provide
 \begin{equation}\label{main7}
 A= \frac{ u (H^2 +\tau^2)}{4 H H_z}.
 \end{equation}
Hence, $A H_z\in \R$. Thus, since by \eqref{main3} we find that $\bar{A} H_z (H+i\tau)\in \R$ (recall that $\kappa-4\tau^2<0$), we can conclude that $H_z^2 (H+i\tau)\in \R$. After writing $z=s+it$, this equation can be rephrased into 
 \begin{equation}\label{main8}
 \tau (H_s^2 -H_t^2) = 2H H_s H_t.
 \end{equation} 

\vspace*{0,3cm}

If $\tau=0$, and thus $\Ek\equiv \H^2(k)\times \R$, then either $H_s\equiv 0$ or $H_t\equiv 0$, i.e. either $H=H(t)$ or $H=H(s)$. In any of these two cases, by the above formulas, we obtain that the fundamental data $\{\landa,H,u,p,A\}$ of the surface depend only on one of the real variables $s$ and $t$. We will label this variable as $r$. We need to show that the surface is an open piece of one of the surfaces in Example \ref{ex2}.


In order to do so, we will assume, up to dilations, that $\kappa =-1$. From \eqref{main4} we infer the existence of a unique (up to $2k\pi$ addition, $k\in\mathbb{Z}$) smooth function $\alpha=\alpha(r)$  such that
\begin{equation}\label{masmain1} 
u=\cos\alpha \qquad 2H=\sin\alpha. 
\end{equation}
As usual, denote by $h$ the last coordinate of the immersion $\psi:\Sigma\to\Ek\equiv\mathbb{H}^2\times\mathbb{R}$. Then 
$A^2=(h_z)^2=\delta (h'(r))^2/4,$ where $\delta=1$ if $r=s$ or $\delta=-1$ if $r=t$. 

Hence, equations \eqref{main1}, \eqref{main6p5} and \eqref{main7}, together with \eqref{masmain1} give that the fundamental data of the immersion coincide with those defined in \eqref{def3.2}. Moreover, from {\bf (C.2)} we infer that $\alpha$ satisfies the differential equation \eqref{eq:3.2} in Example \ref{ex2} and therefore $\psi$ is an open piece of one of the examples described there.
This proves the second part of the theorem.
 
Now, assume that $\tau \neq 0$. Then, using ${\bf (C.2)}$ as well as \eqref{main4} and \eqref{main6p5} it is obtained $$ \left( \frac{u}{2A}\right)_{\bar{z}}= \frac{H^2 +\tau^2}{2|A|^2 (\kappa - 4\tau^2)} - \frac{\lambda u^2 (H+i\tau)}{4A^2}.$$ So, by \eqref{main7} and \eqref{main5} we infer from this last equation that  
 \begin{equation}\label{main9}
(\log (H^2 + \tau^2))_{z\bar{z}} = \frac{ 8 H^2 |H_z|^2}{(H^2 +\tau^2) (4H^2 +\kappa)} + 
\frac{H_z^2 (H+i\tau)(4H^2+\kappa)}{4|H_z|^2(H^2 +\tau^2)}.
 \end{equation}

Summing up now (\ref{main1}), (\ref{main4}), (\ref{main5}), (\ref{main7}), (\ref{main8}) and (\ref{main9}) we conclude that $\psi$ is the one of the surfaces constructed in Example \ref{ex2}. This ends up the proof of Theorem \ref{main}. \\

{\bf A closing remark :} We close this work by analyzing the compact case. So, let $\psi :\Sigma\flecha \Ek$ be a compact surface with holomorphic Abresch-Rosenberg $Qdz^2$, and non-constant mean curvature. If $\tau =0$, it is obvious that there is some point $z_0\in \Sigma$ with $dh (z_0)=0$, i.e. $u(z_0)= \pm 1$. Consequently, from \eqref{main4}, we must have $H(z_0)=0$. However, this is not possible, since when $\tau =0$ we need $H\neq 0$ at every point in order to have $Q dz^2$ well defined. So, compact surfaces with holomorphic $Qdz^2$ are CMC surfaces when $\tau =0$. The same result holds when $P dz^2$ is considered instead of $Q dz^2$ (see \cite{CaLi}).

Finally, assume that $\tau \neq 0$, and suppose that $\Sigma$ has non-zero Euler characteristic. Then there is some point $z_0\in \Sigma$ such that the tangent vector field $T\in \X(\Sigma)$ given in Definition \ref{fundada} vanishes at $z_0$. This implies by ${\bf (C.4)}$ that $u(z_0)=\pm 1$, which is again impossible by \eqref{main4}. Therefore, we conclude that compact surfaces of non-zero Euler characteristic and holomorphic Abresch-Rosenberg differential in $\Ek$ are CMC surfaces.

\section*{Acknowledgements}
Isabel Fernández was partially supported by MEC-FEDER Grant No. MTM2004-00160. Pablo Mira was partially supported by MEC-FEDER, Grant No. MTM2004-02746.

\bibliographystyle{amsplain}

\begin{thebibliography}{10}
\selectlanguage{english}

\bibitem{Abr} U. Abresch, \emph{Generalized Hopf differentials}, preprint.

\bibitem{AbRo} U. Abresch. H. Rosenberg, A Hopf differential for constant mean curvature surfaces in $\S^2\times \R$ and $\H^2\times \R$, {\it Acta Math.} {\bf 193} (2004), 141--174.

\bibitem{AEG} J.A. Aledo, J.M. Espinar, J.A. Gálvez, \emph{Complete surfaces of constant curvature in $\H^2\times \R$ and $\S^2\times \R$}, preprint (http://arxiv.org/abs/math.DG/0510321).

\bibitem{BeTa} D.A. Berdinsky, I.A. Taimanov, Surfaces in three-dimensional Lie groups, \emph{Siberian Math. J.} {\bf 46} (2005), 1005--1019.

\bibitem{CaLi} M.P. Cavalcante, J.H.S. Lira, \emph{Examples and structure of CMC surfaces in some Riemannian and Lorentzian homogeneous spaces}, preprint (http://arxiv.org/abs/math.DG/0511530).

\bibitem{Dan} B. Daniel, \emph{Isometric immersions into $3$-dimensional homogeneous manifolds}, preprint (http://arxiv.org/abs/math.DG/0503500).

\bibitem{FeMi} I. Fernández, P. Mira, \emph{Harmonic maps and constant mean curvature surfaces in $\H^2\times \R$,} preprint (http://arxiv.org/abs/math.DG/0507386).

\end{thebibliography}

\end{document}